\documentclass[a4paper, 12pt]{article}
\usepackage{enumerate, theorem}
\usepackage{amsmath, amsfonts, amssymb}
\usepackage[height=22.5cm, width=15cm]{geometry}
\usepackage[all]{xy}
\usepackage{mathrsfs, yfonts}

\DeclareMathOperator{\id}{id}

 \def\C{{\mathbb C}}

\newcommand{\parag}[1]{\paragraph{\sc{#1.}} }

\newtheorem{thm}{Theorem}[subsection]
\newtheorem{defn}[thm]{Definition}

\newtheorem{prop}[thm]{Proposition}
\newtheorem{lemma}[thm]{Lemma}

\begin{document}

\date{5/04/14}

\author{Daniel Barlet\footnote{Barlet Daniel, Institut Elie Cartan UMR 7502  \newline
Universit\'e de Lorraine, CNRS, INRIA  et  Institut Universitaire de France, \newline
BP 239 - F - 54506 Vandoeuvre-l\`es-Nancy Cedex.France. \newline
e-mail : daniel.barlet@univ-lorraine.fr}.}

\title{ Two  semi-continuity results for the algebraic dimension of  compact complex manifolds.}

\maketitle

\tableofcontents

\parag{Abstract}
Using some relative codimension 1 cycle-space method,  we give, following the ideas of D. Popovici [P.13],  semicontinuity results for the algebraic dimension in a family a compact complex manifolds parametrized by a disc.\\

AMS Classification 2010 : 32G05-32A20-32J10 

\newpage

\section{Statement of results.}

This Note has been inspired by Dan Popovici article [P.13], and we show here that using relative codimension 1 cycle-space, we may give another proof of its theorem 1.2 \ and 1.4  and also
obtain from this proof  much more general results. In fact, very few is new in this approach, because we simply uses the ideas in [P.13] and combined them with the classical tools introduced in [B.75], [B.78] and [C.80]; see also [B-M.14]. This gives a more geometric view on the use of the existence of a strongly Gauduchon form (see definition above) and of the existence of a relative Gauduchon metric 
on a family of compact complex manifolds using a  geometric approach to (relative) algebraic reduction in term of codimension 1  (relative) cycles.\\

We shall prove the following generalizations of the theorem 1.2 and 1.4  of [P.13].

\begin{thm}\label{premier}
Let \ $\pi : \mathcal{X} \to D$ \ be a proper  holomorphic submersion  on an open disc \ $D \in \C$ \ with center \ $0$, where \ $\mathcal{X}$ \ is a connected complex manifold. Assume that for each \ $t \in D$ \ $X_{t} : = \pi^{-1}(t)$ \ is a compact connected complex manifold of dimension \ $n$. Assume also that we have  
\begin{enumerate}[i)]
\item For each \ $t \in D\setminus \{0\}$ \  the manifold \ $X_{t}$ \ is a K{\"a}hler manifold.
\item  For each \ $t \in D\setminus \{0\}$ \ the algebraic dimension of \ $ X_{t}$ \ is at least equal to \ $a \geq 1$.
\item The dimension \ $h^{0,1}(t) : = \dim H^{1}(X_{t}, \mathcal{O})$ \ is independant of \ $t \in D$.
\end{enumerate}
Then the algebraic dimension of the manifold \ $X_{0}$ \ is at least equal to \ $a$.\\
\end{thm}

\begin{thm}\label{second}
Let \ $\pi : \mathcal{X} \to D$ \ be a proper holomorphic submersion on an open disc \ $D \in \C$ \ with center \ $0$, where \ $\mathcal{X}$ \ is a complex manifold. Assume that, for each \ $t \in D^{*}$, the compact connected complex manifold of dimension \ $n$ \  $X_{t} : = \pi^{-1}(t)$ \ is K{\"a}hler. Assume also that \ $X_{0}$ \ admits a smooth \ $2n-2$ \ real \ $d-$closed form such that its \ $(n-1,n-1)$ \ part is positive definite. Then, if for each \ $t \in D\setminus \{0\}$ \ the algebraic dimension of \ $X_{t}$ \ is at least \ $a \geq 1$, the algebraic dimension of \ $X_{0}$ \ is also at least \ $a$.
\end{thm}

\parag{some remarks}
\begin{enumerate}
\item When the map \ $\pi$ \ is k{\"a}hlerian in the sense of [C.81], which implies that \ $X_{t}$ \ is in the class \ $\mathscr{C}$ \ of A. Fujiki\footnote{Note that it  is proved in [V.89] that a compact complex manifold is in the class \ $\mathscr{C}$ \ if and only if it is bimeromorphic to a compact K{\"a}hler manifold.} for all \ $t \in D$, then the semi-continuity result of the algebraic dimension of the fibers of \ $\pi$ \ is proved in [C.81]. \\
But the k{\"a}hlerian assumption for \ $\pi$ \ is rather strong as its implies that \ $X_{0}$ \ is in the class \ $\mathscr{C}$ \ and also properness of the irreducible components of the relative cycles spaces for any dimension of the cycles. Here we only use properness for relative codimension 1 cycles and no K{\"a}hler assumption on \ $X_{0}$. Note that in the absolute case, this compactness for irreducible components of the codimension 1 cycle space is always true (for any compact complex space, see [C.82]) but this is far to be true for smaller cycles of positive dimension.
\item In both cases the key point is to produce a \ $(n-1,n-1)$ \ smooth positive \  $D-$relative form in \ $X$ \ such the integral of it, on the generic member of an analytic family of  relative \ $(n-1)-$cycles giving for generic \ $t \in D\setminus \{0\}$ \ an algebraic reduction for \ $X_{t}$, is locally bounded near \ $0$. In the second case this is an easy consequence of the \ $d-$closedness of this form ; in the first case, it uses only the \ $\partial\bar\partial-$closedness of such a form combined with a nice argument of D. Popovici [P.13] using the constancy of the numbers \ $h^{0,1}(t)$.
\end{enumerate}

\section{Algebraic reduction.}

\subsection{Absolute case.}

For a compact irreducible complex space \ $X$ \ of dimension \ $n$ \ the irreducible components of the complex space \ $\mathcal{C}_{n-1}(X)$ \ are compact  Moishezon spaces. As it is difficult to find an explicit proof of this resul, despite the fact that it is well known for the specialists since 30 years  (see [C.82]), I shall give a proof in the appendix.\\
Note that in codimension \ $> 1$ \ this result is not true in general : see, for instance [C.91].

For any compact irreducible analytic subset \ $\Gamma$ \ of dimension \ $d \geq 1$ \ in \ $\mathcal{C}_{n-1}(X)$ \ there exists a ``natural'' meromorphic map
$$K_{\Gamma}  : X \dashrightarrow \mathcal{C}_{d-1}(\Gamma) $$
called the {\bf Kodaira map} of \ $\Gamma$. This means that there exists a (natural) proper (thanks to [B.78])  modification \ $\tau_{\Gamma} : X_{\Gamma} \to X$ \ and a holomorphic map (also denote \ $K_{\Gamma}$)
$$K_{\Gamma}  : X_{\Gamma} \rightarrow \mathcal{C}_{d-1}(\Gamma)  $$
associating to the generic point \ $x$ \ of \ $X$ \ the \ $(d-1)-$cycle in \ $\Gamma$ \ which is the subset of \ $\Gamma$ \ parametrizing the \ $(n-1)-$cycles containing \ $x$.\\
 We shall denote its image \ $Q_{\Gamma}$. Note that \ $Q_{\Gamma}$ \ is always a compact irreducible Moishezon space. We shall recall some more details on this construction in the relative case in the next paragraph. \\

The algebraic dimension \ $a(X)$ \ of \ $X$ \  is the maximal of the dimension of the \ $Q_{\Gamma}$ \ when \ $\Gamma$ \ is an irreducible component of \ $\mathcal{C}_{n-1}(X)$. In fact, for a given \ $X$ \  there exists a \ $\Gamma$ \ such its Kodaira map is an algebraic reduction (see [C.80] or  [C.81]).

\parag{Claim} In any case, the algebraic dimension of \ $X$ \ is  bigger or equal to the dimension of  \ $Q_{\Gamma}$.
\parag{proof} If the generic member of \ $\Gamma$ is not irreducible (but reduced), there exists a proper finite map\ $ g : \tilde{\Gamma} \to \Gamma$ \ and an analytic family of \ $(n-1)-$cycles parametrized by \ $\tilde{\Gamma}$, with generic irreducible member, obtained by Stein factorization of the projection \ $G_{\Gamma} \to \Gamma$. If \ $\Gamma'$ \ is the image of \ $\tilde{\Gamma}$ \ by the corresponding classifying map, we have a direct image \ $g_{*} : \mathcal{C}_{d-1}(\tilde{\Gamma}) \to \mathcal{C}_{d-1}(\Gamma)$ \ which is proper an finite and induces a genrically finite surjective meromorphic map
$$Q_{\tilde{\Gamma}} \dashrightarrow Q_{\Gamma}.$$
So the dimension of \ $Q_{\tilde{\Gamma}}$ \ is again bounded by the algebraic dimension of \ $X$.\\
In the case of non reduced generic member in the family parametrized by \ $\Gamma$ \ we may consider the pull back of the set theoretic graph \ $\vert G_{\Gamma}\vert$ \ on the normalisation \ $\Gamma_{1}$ \  of \ $\Gamma$ \ which is the graph of an analytic family of \ $(n-1)-$cycles (see [B.75] or [B-M.14]) and replace \ $\Gamma$ \ by the image \ $\Gamma'$ \ of \ $\Gamma_{1}$ \ in \ $\mathcal{C}_{n-1}(X)$ \ by the corresponding classifying map. Then it is easy to see that \ $Q_{\Gamma}$ \ and \ $Q_{\Gamma'}$ \ are bimeromorphic, so again \ $\dim Q_{\Gamma}$ \ is at most  the algebraic dimension of \ $X$.$\hfill \blacksquare$

\subsection{The relative case.}

We consider now a proper surjective holomorphic map \ $\pi : \mathcal{X} \to D$ \ of an irreducible complex space \ $\mathcal{X}$ \ to an open disc \ $D \in \C$ \ with center \ $0$. Let \ $\dim \mathcal{X} = n+1$ ; note that for each \ $t \in D$ \ the fiber \ $X_{t} : = \pi^{-1}(t)$ \ is a compact  complex space with pure dimension\ $n$. It will be convenient to assume that for generic \ $t$ \ in \ $D$ \ the compact space \ $X_{t}$ \ is irreducible\footnote{Using a Stein reduction of \ $\pi$, this assumption is not restrictive.}.\\
Consider now an irreducible \ $D-$proper analytic subset  of dimension \ $d$, $\Gamma \subset \mathcal{C}_{n-1}(\pi)$, of the space of \ $D-$relative \ $(n-1)-$cycles in \ $X$. We say that \ $\Gamma$ \ is a {\bf good filling} \ for \ $\pi : \mathcal{X} \to D$ \ if the following conditions are satisfied
\begin{enumerate}[i)]
\item The generic member of the universal family parametrized by \ $\Gamma$ \ is irreducible.
\item The graph \ $G_{\Gamma} \subset \Gamma \times_{D} \mathcal{X}$ \ of the universal family parametrized by \ $\Gamma$ \ projects surjectively on \ $\mathcal{X}$ \ by the second projection.
\end{enumerate}
Then it is clear that \ $G_{\Gamma}$ \ is irreducible, proper on \ $D$ \ and has dimension \ $d + n-1$. Note that the second condition implies that \ $d + n-1 \geq n+1$, so that \ $d \geq 2$. 

\parag{Claim} The closed analytic subset \ $Y_{\Gamma} \subset \mathcal{X}$ \ which is the locus of the fibers of dimension \ $\geq d-1$ \ of the projection \ $pr : G_{\Gamma} \to \mathcal{X}$ \ has empty interior in each \ $X_{t}, t \in D$.
\parag{proof} If \ $X_{t}$ \ has an irreducible component \ $Z$ \ in \ $Y_{\Gamma}$, then \ $pr^{-1}(Z)$ \ has dimension \ $\geq d + n -1 = \dim G_{\Gamma}$. So, as \ $G_{\Gamma}$ \ is irreducible, we see that this contradicts the fact that \ $pr : G_{\Gamma} \to \mathcal{X}$ \ is surjective.$\hfill \blacksquare$\\

So we have a proper $D-$modification\footnote{This means that it induces a proper modification on {\bf each} \ $X_{t}, t \in D$.} $\tau_{\Gamma} : \mathcal{X}_{\Gamma} \to \mathcal{X}$ \ and a holomorphic \ $D-$map
$$ K_{\Gamma/D} : \mathcal{X}_{\Gamma} \rightarrow \mathcal{C}_{d-2}(\Gamma/D) $$
obtained by composition of the ``fiber-map'' of \ $pr : G_{\Gamma} \to \mathcal{X}$ \ and the $D-$relative direct-image map of  relative\ $(d-2)-$cycles via the $D-$map \ $p : G_{\Gamma} \to \Gamma$ \ (see [B.75] or [B-M.14].

\parag{Claim} The $D-$map \ $K_{\Gamma/D}$ \ is proper, and so \ $Q_{\Gamma}$ \ is proper over \ $D$.
\parag{proof} Let \ $K$ \ be a compact set in \ $\mathcal{C}_{d-2}(\Gamma/D) $ \ and let \ $L$ \ its projection on \ $D$. Then any relative cycle in \ $K$ \ is contained in \ $p^{-1}(L)$ \ where \ $p : \Gamma \to D$ \ is the projection which is assumed to be proper. Then the pull-back of the compact set \ $p^{-1}(L)$ \ on \ $G_{\Gamma}$ \ is a compact set  \ $M$ \ and also its image \ $pr(M)$. Now \ $K_{\Gamma/D}^{-1}(K) \subset pr(M)$ \ and is compact. The propernes of \ $Q_{\Gamma}$ \ on \ $D$ \ is then easy.$\hfill \blacksquare$\\

Denote \ $Q_{\Gamma}$ \ the image of \ $K_{\Gamma/D}$. As for each \ $t \in D$ \ the family \ $\Gamma(t)$ \ is a finite union of irreducible compact analytic subsets in \ $\mathcal{C}_{n-1}(X_{t})$, with at least one good filling set, the compact analytic set \ $Q_{\Gamma}(t)$ \ is a positive dimensional compact Moishezon space.\\

Assume now that for\ $t \in T$, where \ $T$ \ is an uncountable subset of \ $D$, the family \ $\Gamma(t)$ \ gives an algebraic reduction of \ $X_{t}$. Then we have, for such \ $t \in T$ 
$$ \dim Q_{\Gamma}(t) = a(X_{t}) .$$
This implies, by semi-continuity of the dimension of the fibers of the projection \ $Q_{\Gamma} \to D$, that for {\bf each} \ $t$ \ in \ $D$, using the result at the end of paragraph 2.1,  we have
$$ \dim a(X_{t}) \geq \dim Q_{\Gamma}(t) \geq \inf_{t \in T} \dim Q_{\Gamma}(t) = \inf_{t \in T} \  a(X_{t}) .$$

\parag{Conclusion}  Assume that we know the following :
\begin{enumerate}[i)]
\item The irreducible component \ $\Gamma$ \  of \ $\mathcal{C}_{n-1}(\pi)$ \ is proper on \ $D$.
\item For \ $t$ \  in \ $T$ \ where \ $T$ \ is uncountable, we have \ $\dim Q_{\Gamma}(t) =  a$.
\end{enumerate}
Then we conclude that for {\bf each} \ $t \in D$ \ we have \ $a(X_{t}) \geq a $.

\section{A sufficient condition for properness of the irreducible components of the space of relative divisors.}

We begin by a definition which is equivalent to the notion of a ``strongly Gauduchon metric'' introduced in [P.13].

\begin{defn}
Let \ $X$ \ be a reduced complex space of pure dimension \ $n$. A  \ $2n -2$ \ smooth form \ $\omega$ \ on \ $X$ \ will be called a  sG-form on \ $X$ \ when it satisfies:
\begin{enumerate}[i)]
\item The form \ $\omega$ \ is \ $d-$closed.
\item The \ $(n-1,n-1)$ \ part of \ $\omega$ \ is positive definite on \ $X$.
\end{enumerate}
\end{defn}

Recall that the strict  positivity above means that in a local embedding in an open set \ $U$ \ of \ some \ $\C^{N}$, the \ $(n-1,n-1)$ \ part of \ $\omega$ \ may be induced on \ $X$ \ by a smooth, strictly positive in the sense of Lelong \ $(n-1,n-1)$ \ form on \ $U$.\\
We also need the relative version of this notion :

\begin{defn}
Let \ $ \pi :\mathcal{X} \to S$ \ be a surjective proper $n-$equidimensional morphism of reduced complex spaces. A smooth $S-$relative \ $2n-2$ \ form \ $\omega_{/S}$ \ will be called a $S-$relative sG-form if it induces a sG-form on each \ $X_{s} : = \pi^{-1}(s)$.
\end{defn}

Note that in the absolute case the irreducible components of the space \ $\mathcal{C}_{n-1}(X)$ \ of compact \ $(n-1)-$cycles of a compact pure \ $n-$dimensional complex space \ $X$ \ are always compacts, thanks to the existence of a meromorphic algebraic reduction \ $X \dashrightarrow Y$, where \ $Y$ \ is a projective variety (see [C.81] \ and [C.82]). The existence of a $S-$relative sG-form for a proper $n-$equidimensional surjective map \ $\pi : \mathcal{X} \to S$ \ implies a relative version of this result.

\begin{prop}\label{propre}
Let \ $ \pi :\mathcal{X} \to S$ \ be a surjective proper $n-$equidimensional morphism of reduced complex spaces admitting a $S-$relative sG-form \ $\omega$. Then each connected component of the \ $S-$relative cycle space \ $\mathcal{C}_{n-1}(\pi)$ \ is proper over \ $S$.
\end{prop}

\parag{Proof} For a \ $(n-1)-$relative cycle  \ $C$ \ of \ $\pi$ \ define \ $F(C) : = \int_{C} \omega$. This is a continuous function of \ $C$ \ (see [B.75] or [B-M.14]) and  on any given compact set in \ $S$, this function is bigger than \ $\epsilon$\ times the volume of \ $C$ \ for a continuous hermitian metric defined on \ $\mathcal{X}$, because in the integration, only the \ $(n-1,n-1)$ \ part of \ $\omega_{/S}$ is relevant. As the function \ $F$ \ is locally constant on the fibers of the projection \ $\mathcal{C}_{n-1}(\pi) \to S$ \ thanks to the \ $d_{/S}-$closeness of \ $\omega_{/S}$, this is enough to prove the properness of the projection \ $\mathcal{C}_{n-1}(\pi) \to S$ \ thanks to E. Bishop's theorem (for instance, see [B-M.14]). $\hfill \blacksquare$

\begin{lemma}\label{Stab.}{\rm [see [P.13]]}
Let \ $\pi : \mathcal{X} \to D$ \ be a proper holomorphic family of compact connected complex manifolds of dimension \ $n$ \ parametrized by an open disc with center \ $0$ \ in \ $\C$. Assume that \ $X_{0} : = \pi^{-1}(0)$ \ has a G-form \ $\omega_{0}$. Then we may find a small open disc \ $D' \subset D$ \ with center \ $0$ \ and a relative G-form \ $\omega$ \ on \ $\pi^{-1}(D')$ \ inducing \ $\omega_{0}$ \ on \ $X_{0}$.
\end{lemma}

\parag{Proof} Thanks to Ehresmann theorem, there exists an open disc \ $D_{1}\subset D$ \ with center \ $0$ \ and a \ $\mathscr{C}^{\infty}$ \ trivialisation \ $\theta : \pi^{-1}(D_{1}) \to D_{1}\times X_{0}$ \ of the fibration \ $\pi$ \ inducing the identity on \ $X_{0}$. Define \ $\omega : = \theta^{*}(\omega_{0})$. As \ $\omega_{0}$ \ is \ $d-$closed, so is \ $\omega$. We shall consider \ $\omega$ \ as a relative \ $d-$closed form. It induces \ $\omega_{0}$ \ on \ $X_{0}$. As the complex structure of \ $X_{t}$ \ varies continuously with \ $t \in D_{1}$, the \ $(n-1,n-1)$ \ part of the relative form \ $\omega$ \ varies continuously. As it is definite positive at \ $t = 0$ \ there exists an open disc \ $D' \subset D_{1}$ \ with center \ $0$ \ where it stays definite positive on each fiber ; so \ $\omega$ \ induces a relative G-form on \ $\pi^{-1}(D')$. $\hfill \blacksquare$.\\

In the situation of this corollary define \ $\mathcal{X}' : = \pi^{-1}(D')$. Now the proposition above gives that the connected components of \ $\mathcal{C}_{n-1}(\mathcal{X}')$ \ are proper over \ $D'$.

\parag{Remark} In the situation of theorems \ref{premier} and \ref{second} we may apply this lemma to see that any irreducible component of \ $\mathcal{C}_{n-1}(\pi^{*})$ \ is proper over \ $D^{*}$, where \ $D^{*} = D \setminus \{0\}$ \ and \ $\pi^{*}$ \ is the restriction to \ $\pi^{-1}(D^{*})$, because we assume that each \ $X_{t}, t \in D^{*}$, is K{\"a}hler, so admits a G-form.

\section{Adaptation of two key results of [P.13].}

We shall adapt in this section two key results in [P.13] in order to use them in our proofs. 

\subsection{Adaptation of the remark 2.1}

Let \ $\pi : \mathcal{X} \to D$ \ be a proper holomorphic submersion of a \ $(n+1)-$dimensional connected complex manifold \ $\mathcal{X}$ \ to a disc. To follow the point of view in [P.13], we trivialize \ $\pi$ \ as a 
 \ $\mathscr{C}^{\infty}-$submersion\footnote{this is not restrictive, up to localize the situation on \ $D$; see the begining of the proof of the lemma \ref{Stab.}.} and we shall assume that \ $\mathcal{X} \equiv X \times D$, $\pi$ given by the second projection, the complex manifold  \ $X_{t} : = \pi^{-1}(t)$ \ is then defined by a complex structure \ $J_{t}$ \ on \ $X$ \ which depends smoothly on \ $t \in D$.\\
 We shall denote by \ $\mathcal{X}^{*} : = \pi^{-1}(D^{*})$ \ where \ $D^{*} : = D \setminus \{0\}$, and \ $\pi^{*} : \mathcal{X}^{*}  \to D^{*}$ \ the restriction of \ $\pi$. Note that this morphism is (locally on \ $D^{*}$)  a K{\"a}hler morphism, so that the irreducible components of \ $\mathcal{C}_{n-1}(\pi^{*})$ \ are proper on \ $D^{*}$ \ (see the final remark of the previous paragraph).

\begin{prop}\label{choice}
In the situation above, assume that the algebraic dimension \ $a(X_{t})$ \ of \ $X_{t}$ \ is at least equal to \ $a \geq 1$ \ for all \ $t \in D\setminus \{0\}$. Then there exists an irreducible component \ $\Gamma^{*}$ \ of the space \ $\mathcal{C}_{n-1}(\pi^{*})$ \ with the following properties :
\begin{enumerate}[i)]
\item \ $\Gamma^{*}$ \ is proper and surjective over \ $D^{*}$.
\item \ $\Gamma^{*}$ \ fills up \ $\mathcal{X}^{*} $ \ (i.e. the graph \ $G_{\Gamma^{*}}$ \.of \ $\Gamma^{*}$ \ surjects onto \ $\mathcal{X}^{*}$).
\item The generic cycle in \ $\Gamma^{*}$ \ is irreducible.
\item The dimension of \ $Q_{\Gamma^{*}}$ \ is at least \ $a + 1$ \ and the projection \ $Q_{\Gamma^{*}} \to D^{*}$ \ is proper and surjective.
\end{enumerate}
\end{prop}

\parag{proof} Note first that, as any irreducible component of \ $\mathcal{C}_{n-1}(\pi^{*})$ \ is proper over \ $D^{*}$, the only condition in  i) is the surjectivity.\\
For each \ $t \in D^{*}$ \ there exists an irreducible \ $(n-1)-$cycle \ $C_{t}$ \ in \ $X_{t}$ \ which is the generic member of a compact  irreducible analytic subset \ $\Gamma(t)$ \  in \ $\mathcal{C}_{n-1}(X_{t})$ \ such that the image \ $Q_{\Gamma}(t)$ \ of the corresponding Kodaira map as dimension at least \ $a$. Then \ $\Gamma(t)$ \ is a compact analytic subset of some irreducible component \ $\Gamma^{*}$ \ of \ $\mathcal{C}_{n-1}(\pi^{*})$. As there are only countably many \ $\Gamma^{*}$ \ and uncountably many \ $t \in D^{*}$, there exists at least one \ $\Gamma^{*}$ \ which contains uncountably many such \ $\Gamma(t)$. Then this \ $\Gamma^{*}$ \ must satisfies i). The corresponding \ $Q_{\Gamma^{*}}$ \ has uncountably many fibers over \ $D^{*}$ \ which have dimension at least \ $a$. So the condition iv) is also fullfilled, as properness  of \ $Q_{\Gamma^{*}}$ \ on \ $D^{*}$ \ is automatic form the properness of \ $\Gamma^{*}$ \ on \ $D^{*}$. This implies condition ii) because the fact that \ $Q_{\Gamma^{*}}$ \ surjects to \  $D^{*}$ \ implies that the generic fiber of \ $\Gamma^{*}$ \ over \ $D^{*}$ \ has dimension \ $\geq 1$. The condition iii) is also clear.$\hfill \blacksquare$
\bigskip

 In this situation we shall denote by \ $\alpha(\Gamma^{*}) \in H^{2}(X, \mathbb{Z})$, the (topological) fundamental class of any (relative) cycle in the family parametrized by \ $\Gamma^{*}$.

\subsection{Adaptation of the proposition 3.1.}

Recall first that in our situation there exists, thanks to [G.77],  a smooth relative Gauduchon metric on \ $\mathcal{X}$. For our purpose we shall simply use the fact that there exists a smooth family \ $\gamma_{t}$ \ of positive definite \ $(1,1)-$forms on \ $X_{t}, t \in D$ \ with the condition that \ $\partial_{t}\bar\partial_{t} \ \gamma_{t}^{\wedge (n-1)} = 0$ \ for each \ $t\in D$. Then the variant of the proposition 3.1. of [P.13] we shall use is the following.
\begin{prop}\label{contr. vol.}
Let \ $\alpha \in H^{2}(X, \mathbb{Z})$ \ be the fundamental class of the relative \ $(n-1)-$cycles in the family parametrized by \ $\Gamma^{*}$ \ chosen as in the proposition \ref{choice}. Note \ $[C_{t}]$ \  the integration current on \ $X_{t}$ \ of a member of the family \ $\Gamma^{*}$ \ contained in \ $X_{t}$ \ for some \ $t \in D^{*}$. If \ $h^{0,1}(t)$ \ is independent of \ $t \in D$, there exists a constant \ $C > 0$ \ such that for \ $t \in D^{*}$ \ near enough \ $0$ \ we have
$$ < [C_{t}], \gamma_{t}^{\wedge n-1} > \  \leq C < +\infty $$
where \ $\gamma_{t}^{\wedge(n-1)}$ \ is the \ $(n-1)-$th exterior power of a relative Gauduchon metric on \ $\mathcal{X}$.
\end{prop}

\parag{proof} Let \ $\tilde{\omega}$ \ be a smooth \ $d-$closed real \ $2-$form on \ $X$ \ which is a de Rham representative of  the class \ $\alpha \in H^{2}(X, \mathbb{Z})$. So for each \ $t \in D^{*}$ \ there exists a real \ $1-$current \ $\beta_{t}$ \ on \ $X_{t}$ \ such that \ $[C_{t}] = \tilde{\omega} + d\beta_{t}$. As \ $\beta_{t}$ \ is real, we have \ $\beta_{t} = \beta_{t}^{1,0} + \overline{\beta_{t}^{1,0}}$.\\
Note that type consideration shows that \ $\beta_{t}^{0,1}: = \overline{\beta^{1,0}}$ \ is solution of the equation
$$ \bar\partial_{t} \beta_{t}^{0,1} = - \tilde{\omega}_{t}^{0,2} .$$
So following [P.13], we shall define \ $\tilde{\beta}^{0,1}$ \ as the (unique) solution of this equation with minimal \ $L^{2}$ \ norm defined by the Gauduchon metric \ $\gamma_{t}$ \ on \ $X_{t}$. Defining \ $\tilde{\beta}_{t} = \overline{\tilde{\beta}^{0,1}} + \tilde{\beta}^{0,1} $ \ we have now that the real  current
$$  [C_{t}] - \tilde{\omega}_{t} - d(\tilde{\beta}_{t}) $$
is \ $d-$exact and of type \ $(1,1)$ \ on \ $X_{t}$. Then, as the complex compact manifold \ $X_{t}$ \ is K{\"a}hler, it satisfies the \ $\partial\bar\partial-$lemma (see for instance [V.86]), and there exists a \ $0-$current \ $\varphi_{t}$ \ on \ $X_{t}$ \ such that
$$  [C_{t}]  =  \tilde{\omega}_{t} + d(\tilde{\beta}_{t}) + i\partial_{t}\bar\partial_{t} \varphi_{t} .$$
An easy consequence is now, as \ $\partial_{t}\bar\partial_{t} \ \gamma_{t}^{\wedge (n-1)} = 0$, that
$$  < [C_{t}], \gamma_{t}^{\wedge n-1} > = <  \tilde{\omega}_{t} , \gamma_{t}^{\wedge n-1} > + < d(\tilde{\beta}_{t}),  \gamma_{t}^{\wedge n-1} > .$$
As the function \ $t \mapsto  <  \tilde{\omega}_{t} , \gamma_{t}^{\wedge n-1} > $ \ is continuous on \ $D$, to bound the left handside, it is enough to bound near \ $t = 0$ \ the term
$$  < d(\tilde{\beta}_{t}),  \gamma_{t}^{\wedge n-1} >  = \pm < \tilde{\beta}_{t}, d\gamma_{t}^{\wedge n-1} >  .$$
Then it is enough to follow the argument concluding the proof of the proposition 3.1 in  [P.13], showing that under the assumption that \ $h^{0,1}(t)$ \ is constant, $t \mapsto \tilde{\beta}_{t}$ \ depends continuously of \ $t \in D$ \ to conclude. $\hfill \blacksquare$\\

An obvious consequence of this result is that the volume (for the relative Gauduchon metric we have fixed, but this property is independant of the choice of a continuous relative hermitian metric) of members in \ $\Gamma^{*}$ \ is uniformely bounded, up to shrink the disc \ $D$ \ around \ $0$. This implies that the closure \ $\Gamma$ \ of \ $\Gamma^{*}$ \ in \ $\mathcal{C}_{n-1}(\pi)$ \ which is an irreducible analytic subset, is proper over \ $D$ \ thanks to Bishop's theorem (see [B-M.14]). \\

\parag{Proofs of the two theorems}
In both cases we have an irreducible analytic subset \ $\Gamma$ \ of \ $\mathcal{C}_{n-1}(\pi)$ \ which is proper surjective over \ $D$, with an irreducible generic member, such that its graph surjects to \ $\mathcal{X}$ \ and such that the corresponding \ $Q_{\Gamma}$ \ has generic fibers of dimension \ $\geq a$ \ over \ $D$. Then its fiber at \ $t = 0$ \ has dimension at least \ $a$ \ which implies that the algebraic dimension of \ $X_{0}$ \ is at least \ $a$, thanks to the conclusion of the section 2.$\hfill \blacksquare$

\newpage

\section{Appendix}

We give here a proof of the following  (classical) statement :

\begin{thm}
Let \ $X$ \ be a compact irreducible complex space of dimension \ $n$. Then the irreducible components of the space \ $\mathcal{C}_{n-1}(X)$ \  of \ $(n-1)-$cycles in \ $X$ \ are compact and Moishezon.
\end{thm}

\parag{proof} Using Hironaka's desingularization theorem it is not restrictive to assume that \ $X$ \ is a compact connected manifold of dimension \ $n$ \ and that we have a holomorphic surjective map \ $r : X \to P$ \ where \ $P$ \ is a projective manifold of dimension \ $a$ \ where \ $a$ \ is the algebraic dimension of \ $X$. This comes from the fact that by a proper modification of \ $X$ \ we add only finitely many effective irreducible \ $(n-1)-$cycles, and then, for any irreducible component \ $\Gamma$ \ of the space of \ $(n-1)-$cycle of our initial \ $X$, there exists an irreducible component \ $\tilde{\Gamma}$ \ of  space of \ $(n-1)-$cycles of the smooth modification of \ $X$ \ such that the direct image of cycles gives a modification \ $\tilde{\Gamma} \to \Gamma $. \\
Recall now that the number of non polar effective irreducible divisors in \ $X$ \ is bounded thanks to [C.82]. Let now \ $\Gamma$ \ be an irreducible component of \ $\mathcal{C}_{n-1}(X)$ \ of positive dimension with generic member irreducible\footnote{It is enough to treat this case for our result.}. Then each member of the corresponding family is polar and so has a \ $(a-1)-$dimensional image in \ $P$. This means that the image of the graph \ $G_{\Gamma} \subset \Gamma \times X$ \ by the map \ $\id_{\Gamma}\times r$ \ is proper and \ $(a-1)-$equidimensional on \ $\Gamma$. Up to replace \ $\Gamma$ \ by its normalization, we have the graph of an analytic family of \ $(a-1)-$cycles in \ $P$. As \ $\Gamma$ \ is irreducible and \ $P$ \ projective, the volume of these cycles in \ $P$ \ is uniformely bounded.\\
Now using any continuous hermitian metric on \ $X$, we have an uniform bound for the volume of the generic fibers of \ $r$ \ thanks to [B.78]. Then a Fubini type argument implies that the volume of the generic member of our initial family, which is bounded by the volume of the pull-back by \ $r$ \ of its image by \ $r$ \ is uniformely bounded. Now Bishop's theorem implies that \ $\Gamma$ \ is relatively compact in \ $\mathcal{C}_{n-1}(X)$, as it is closed, it is compact. To conclude we have to remember that the normalization of \ $\Gamma$ \ dominates a compact analytic subspace of an irreducible component of \ $\mathcal{C}_{a-1}(P)$ \ and that this map is generically finite, because the pull-back by \ $r$ \ of an irreducible effective divisor in \ $P$ \ contains only finitely many irreducible effective divisors in \ $X$. So \ $\Gamma$ \ is a compact irreducible Moishezon space.$\hfill \blacksquare$

\newpage

\section{References}

\begin{enumerate}

\item[[B.75]] Barlet, D. {\it Espace analytique r\'{e}duit des cycles analytiques complexes compacts d'un espace analytique complexe de dimension finie}, Fonctions de plusieurs variables complexes, II (S\'{e}m. Franois Norguet, 1974-1975), pp.1-158. Lecture Notes in Math., Vol. 482, Springer, Berlin, 1975

\item[[B.78]] Barlet, D. {\it Majoration du volume des fibres g\'{e}n\'{e}riques et forme g\'{e}om\'{e}trique du th\'{e}or\`{e}me d'aplatissement}, S\'{e}minaire Pierre Lelong-Henri Skoda (Analyse) 1978/79  pp.1-17, Lecture Notes in Math., 822, Springer, Berlin-New York, 1980.

\item[[B-M.14]] Barlet,D. et Magnusson,J. {\it Cycles analytiques complexes I : Th\'{e}or\`{e}mes de pr\'{e}paration des cycles} \`{a} para\^{i}tre dans la s\'{e}rie ``Cours Sp\'{e}cialis\'{e}s'' SMF 2014.

\item[[C.80]] Campana,F. {\it Application de l'espace des cycles \`{a} la classification bim\'{e}romorphe des espaces analytiques k{\"a}hl\'{e}riens compacts}, Revue de l'Inst. E. Cartan \ $n^{0} 2$ (Nancy), pp.1-162.

\item[[C.81]] Campana, F. {\it R\'{e}duction alg\'{e}brique d'un morphisme faiblement K{\"a}hl\'{e}rien propre et applications}, Math. Ann. 256 (1981), no. 2, pp.157-189.

\item[[C.82]] Campana,F. {\it Sur les diviseurs non-polaires d'un espace analytique compact}, J. Reine Angew. Math. 332 (1982), pp.126-133.

\item[[C.91]] Campana,F. {\it The class \ $\mathscr{C}$ \ is not stable by small deformations}, Math. Ann. 290 (1991), no. 1, pp.19-30.

\item[[G.77]] Gauduchon, P. {\it Le th\'{e}or\`{e}me de l'excentricit\'{e} nulle}, C.R.Acad.Sci. Paris, S\'{e}r. A 285, pp.387-390 (1977).

\item [[P.13]] Popovici,D.  {\it Deformation limits of projective manifolds: Hodge numbers and strongly Gauduchon metrics}, Invent. Math. 194 (2013), no. 3, pp.515-534. 

\item[[V.86]] Varouchas, J. {\it Propri\'{e}t\'{e}s cohomologiques d'une classe de vari\'{e}t\'{e}s analytiques complexes compactes}, S\'{e}minaire d'analyse P. Lelong-P. Dolbeault-H. Skoda,  1983/1984, pp.233-243, Lecture Notes in Math., 1198, Springer, Berlin, 1986.

\item[[V.89]] Varouchas, J. {\it K{\"a}hler spaces and proper open morphisms}, Math. Ann. 283 (1989), no. 1, pp.13-52.

\end{enumerate}

\end{document}